\documentclass[11pt]{article}

\input{amssym.def}
\input{amssym}
\usepackage{eufrak}

\newtheorem{theorem}{Theorem}

\newtheorem{lemma}[theorem]{Lemma}
\newtheorem{proposition}[theorem]{Proposition}
\newtheorem{remark}[theorem]{Remark}

\def\Uq{U_q(sl(n))}
\def\uq{U_q(sl(2))}

\def\ot{\otimes}
\def\h{{\hbar}}
\def\K{{\Bbb K}}
\def\CC{{\Bbb C}}
\def\R{{\Bbb R}}
\def\X{{\cal X}}
\def\Y{{\cal Y}}
\def\Z{{\cal Z}}
\def\B{{\cal B}}
\def\C{{\cal C}}
\def\RR{{\textsc R}}
\def\H{{\cal H}}
\def\LL{{\cal L}}
\def\SL{{\cal SL}}

\def\ad{{\rm ad\, }}
\def\End{{\rm End}}
\def\Vect{{\rm Vect}}

\def\vv{V^{\otimes 2}}

\def\De{\Delta}
\def\om{\omega}
\def\Om{\Omega}
\def\M{\cal{M}}
\def\N{\cal{N}}
\def\Ker{\rm Ker}

\def\grad{\rm grad}
\def\div{\rm div}
\def\rot{\rm rot}
\def\al{\alpha}
\def\Mat{\rm Mat}
\def\ah{A_\h}

\def\sr{S_r^2}
\def\Mw{\rm Mw}

\def\grad{\rm grad}
\def\diag{\rm diag}

\def\div{\rm div}
\def\Pois{\rm Pois}
\def\Cas{\rm Cas}

\def\Id{\rm Id}
\def\span{\rm span}

\def\lhq{\ifmmode {\cal L}(q,\hbar)\else ${\cal L}(q,\hbar)$\fi}

\def\lrq{{\cal L}(R_q)}

\def\lqh{\ifmmode {\cal L}(q,\hbar)\else ${\cal L}(q,\hbar)$\fi}
\def\Hr{{\rm H}_r^2}

\def\qq{q^{-1}}

\def\Tr{{\rm Tr}}

\def\al{{\alpha}}

\def\be{\begin{equation}}
\def\ee{\end{equation}}

\begin{document}

\makeatletter
\renewcommand{\theequation}{{\thesection}.{\arabic{equation}}}
\@addtoreset{equation}{section} \makeatother

\title{Maxwell operator on q-Minkowski space\\
 and q-hyperboloid algebras}
\author{Antoine Dutriaux\thanks{dutriaux@univ-valenciennes.fr},
Dimitri Gurevich\thanks{gurevich@univ-valenciennes.fr}\\
{\small\it LAMAV, Universit\'e de Valenciennes, 59304 Valenciennes,
France}}

\maketitle

\begin{abstract}
We introduce  an analog of the Maxwell operator on a q-Minkowski
space algebra (treated as a particular case of the so-called
Reflection Equation Algebra) and on certain of its quotients. We
treat the space of "quantum differential forms" as a projective
module in the spirit of the Serre-Swan approach. Also, we use
"braided tangent vector fields" which are q-analogs of Poisson
vector fields associated to the  Lie bracket $sl(2)$.
\end{abstract}

{\bf AMS Mathematics Subject Classification, 2000:} 17B37, 81R50


\section{Introduction}

The main goal of this paper is to define a q-analog of the Maxwell
operator on some noncommutative (NC) algebras. Namely, we are
dealing with three of them: the q-Minkowski space algebra
$\K_q[\R^4]$, quantum (braided, or q-)hyperboloid algebra
$\K_q[\Hr]$, and an intermediate algebra $\K_q[\R^3]$. The
q-Minkowski space algebra (as defined in \cite{MM}, \cite{M},
\cite{K}) is a particular case of a so-called reflection equation
algebra (REA), the
others are its quotients.   Observe that the REA was introduced  in the
early 90's by S.Majid under the name of braided matrix algebra (cf. the cited papers and the references therein). He
also defined a Hermitian structure in it. Here we do not consider
such a structure\footnote{The problem of defining involution operators  in
"braided algebras" will be discussed in subsequent papers.
(Hereafter the term "braided" stands for the REA and related algebras and objects.)
We would like only to note that the  q-Minkowski space algebra endowed with the mentioned Hermitian structure
cannot be treated as a real vector space. }.

The above algebras are deformations of their classical counterparts
$\K[\R^4]$, $\K[\Hr]$, and $\K[\R^3],$ respectively. Hereafter
$\K=\CC$ (or $\R)$ is the ground field, the notation $\K[\M]$
stands for the coordinate algebra of a given regular affine
algebraic variety $\M$, and \be \Hr=\{(b,\,h,\,c) \in \R^3\,|\,
2bc+{h^2\over 2}=r^2,\, r\not=0 \} \label{hyper} \ee is a
hyperboloid\footnote{If $\K=\R$ and $r$ is real, we get a
one-sheeted hyperboloid. If $r$ is purely imaginary, we get a
two-sheeted hyperboloid. However, if $\K=\CC$ we allow $r$ to take
any non-trivial value. It should be emphasized that we prefer to
deal with a q-analog of a hyperboloid (and not of a sphere -- the
so-called Podles' sphere) since it cannot be realized as a real
algebra.}. Moreover, all deformed algebras we are dealing with,
can be endowed with an action of the quantum group (QG) $\uq$
compatible with their product in the usual way.

By passing from the  classical algebras to their quantum or
braided analogs we want to
simultaneously deform certain differential operators defined on
the initial algebras or some vector bundles over them. Moreover,
if a given operator is covariant w.r.t. to a group $G$ we require
its deformed counterpart to be covariant w.r.t. the corresponding
QG.
 The simplest operator which can be "q-deformed" is the  Laplace (or Laplace-Beltrami) operator
 on one of the mentioned algebras.
 Recall that this operator is associated to a (pseudo-)Riemannian metric on a given regular variety.
Thus, if the metric $g$ which comes in its definition is constant
of the signature $(1|3)$ on the space $\R^4$ the corresponding
Laplace operator (also called d'Alembertian)  is
$\partial_t^2-\partial_x^2-\partial_y^2-\partial_z^2$ where
$(t,\,x,\,y,\,z)$ are Cartesian coordinates in this space. If
${\M}=\Hr$ and $g\in \Om^2(\Hr)$ is an $SL(2)$-invariant metric
the corresponding  Laplace operator equals (up to a factor) to
 the quadratic Casimir operator coming from
the enveloping algebra $U(sl(2))$ whereas the hyperboloid is
treated as an orbit $\Hr \hookrightarrow sl(2)^*$ of the coadjoint
action of the Lie algebra $sl(2)$. Vector fields arising from this
coadjoint action are tangent to all orbits in $sl(2)^*$ and we
call them
 {\em tangent vector fields}. (Also, they are Poisson vector fields w.r.t. the linear Poisson-Lie bracket
 defined on the space $sl(2)^*$.)

Thus, in order to define braided
 analogs of the Laplace  operator on the quantum algebras in question
we should first introduce braided analogs of vector fields. More
precisely, we need analogs of tangent vector fields while dealing
with the algebra $\K_q[\Hr],$ and analogs of partial derivatives
when dealing with the algebras $\K_q[\R^3]$ and $\K_q[\R^4]$.
 The problem of defining a braided analog of the Maxwell operator is even more complicated since
we should first introduce braided analogs of the spaces of
differential forms $\Om^1(\Hr)$,  $\Om^1(\R^3)$, and
$\Om^1(\R^4)$. Recall that the Maxwell operator is defined on a
given variety $\M$ as follows: \be \omega\to
\partial\,d\,\omega,\,\,\, \partial=*^{-1}\,d \, * .\label{Max}\ee
Here $\om\in \Om^1(\M)$ is a one-form on  $\M$, $d$ is the de Rham
operator and $*$ is the Hodge operator. (Note that on the
classical Minkowski  space the Maxwell system is initially defined
on the space $\Om^2(\R^4)$  but it can be easily reduced to the
operator above. Also, note that the conventional definitions of
the Maxwell and Laplace operators
 differ from ours by a sign. We disregard this.)

There are several  approaches to the problem of defining
analogs of differential forms and of
 the de Rham operator on a given noncommutative algebra $A$. The first approach consists of considering
{\em universal differential forms} without any commutation
relations (e.g., $a\,db=db\,a$) between "functions" $a\in A$ and
"differentials" $d b$ but with the preservation of the Leibnitz
rule. The corresponding differential algebra is much bigger than
the algebra of usual differential  forms even if the initial
algebra is the coordinate algebra $A=\K[\M]$ of a regular variety
$\M$.

If an algebra $A$ is related to a braiding (say, it is a so-called
RTT algebra or an REA, see section 4) one looks for an extension
of the braiding coming in the definition of $A$ onto the space of
differential forms. Such an extended braiding  enables one to
relate the elements of the form $a\,db$ and $db\, a,$ and to
reduce the space of universal differential forms  to the
"classical size". In the case of the quantum analog of the group
$GL(n)$ this can even be done with preservation of the Leibnitz
rule (cf. \cite{W, K, IP}) while for the quantum analog of the
group
  $SL(n)$ this rule has to be dropped (cf. \cite{FP}).

The third approach, due to A. Connes, is based on the  notion of
spectral triples. In the framework of this approach the role of
differential forms is played  by the (classes of)  Hochschild
cycles \cite{C}.

Nevertheless, all these approaches do not enable one to define a smoothly deformed space of differential
forms on a quantum hyperboloid algebra
 $\K_q[\Hr]$.
As it was observed in \cite{AG} the space of differential forms
$\Om^1(\Hr)$  can be smoothly deformed
 into a quantum one $\Om^1_q(\Hr)$   as a
one-sided module.  However, opposite  to the classical case
($q=1$), if we convert this one-sided $\K_q[\Hr]$-module into a
two-sided module via an extension of the initial braiding we
reduce the size of the space of differential forms.

Following \cite{AG, A}, we treat the spaces of braided
differential forms as one-sided projective modules\footnote{Note
that in general, due to the Serre-Swan approach, any vector bundle
on  a regular affine variety can be realized as a projective
module. Recall that such (say, right)
  $A$-module
 for a given algebra $A$ is of the form $e\,A^{\oplus n}$
where  $e\in {\Mat}_n(A)$ is  an idempotent and ${\Mat}_n(A)$ stands for the space of
$n\times n$ matrices with entries from $A$. \\
As it was shown in \cite{R}, if $\ah$ is a formal deformation of a
commutative algebra $A,$ any idempotent $e\in {\Mat}_n(A)$ can be
deformed into an idempotent $e_\h\in {\Mat}_n(\ah)$. Thus, we get
a one-sided projective $\ah$-module $e_\h\,\ah^{\oplus n}$ which
is a formal deformation of the initial $A$-module. (Nowadays,
there is a known an explicit formula for such a deformed idempotent,
cf. \cite{BB}.)}. In order to do so, we use the following
remarkable property of the q-Minkowski space algebra $\K_q[\R^4]$.
There exists a series of the Cayley-Hamilton (CH) polynomial
identities for some matrices  with entries from the algebra
$\K_q[\R^4]$. The coefficients of the CH polynomials are central
in this algebra, and they become scalar if we switch to the
quotient $\K_q[\Hr]$. It is  a somewhat standard trick to use
these polynomials   to construct a set of idempotents and
corresponding projective modules. Thus, we can explicitly deform
any $SL(2)$-equivariant vector bundle on a hyperboloid $\Hr,$
realized as a projective module, to its braided counterpart. By
applying this scheme to the space  $\Om^1(\Hr)$ we get its braided
analog $\Om^1_q(\Hr),$ realized as a projective
$\K_q[\Hr]$-module. Moreover, it is $\uq$-equivariant (covariant).

Besides, we define braided analogs of tangent vector fields without any
form of the Leibnitz rule. In order to do that,
we use another remarkable property of the q-Minkowski space
algebra $\K_q[\R^4]$. Let $\LL$ be the space spanned by the
generators of this algebra. There exists a braided analog
$$[\,,\,\,]_q:\LL^{\ot 2}\to\LL$$
 of the $gl(2)$  Lie bracket such that (a slight modification of)
the algebra $\K_q[\R^4]$ can be regarded as the enveloping algebra
of the corresponding q- (or braided) Lie algebra. We refer the
reader to \cite{GPS2} for further explanations (applicable in a
much more general setting). Below, we only need a braided analog
$[\,,\,\,]_q:\SL^{\ot 2}\to\SL$ of the Lie algebra $sl(2)$ where
$\SL$ is a 3 dimensional subspace of the space $\LL$. By using
this q-bracket we define "braided tangent vector fields" following
the classical pattern. Moreover, there is a braided analog of the
quadratic Casimir element on the space $\SL^{\ot 2}$. Representing
it by
 the above "braided tangent vector fields" we get an analog  of the Casimir operator on a q-hyperboloid.

A braided analog of the Maxwell operator on a q-hyperboloid  is
also defined via braided tangent vector fields as a proper
deformation of the Maxwell operator acting on $\Om^1(\Hr)$.
However, as we have said above, the space of braided differential
forms on a q-hyperboloid is treated as a one-sided projective
$\K_q[\Hr]$ module. In order to do this we first apply this scheme
to the Maxwell operator on a usual hyperboloid (which seems to be
new even in this classical setting).

In order to get the Maxwell operator on the algebra $\K_q[\R^3],$
in addition to "braided tangent vector fields", we must also use
 the derivative in $r$ where $r$ is an analog of the radius (it comes
 in a
parametrization of quantum hyperboloids and in a sense belongs to
the algebraic extension of the center of the algebra
$\K_q[\R^3]$). We assume that this derivative  has the classical
properties, in particular, it satisfies the Leibnitz rule.
 Thus, we relate the braided Laplace and Maxwell  operators
on the algebra $\K_q[\Hr]$ and  on the algebra  $\K_q[\R^3]$ in a
way similar to the classical one: the  operators on the former algebra are
restrictions of those on the latter one. Emphasize that the methods of
defining partial derivatives on the q-Minkowski space algebra via
"braided Leibnitz rule"
(based on a transposition of "functions" and "partial derivatives" as in \cite{K},  \cite{IP}, \cite{FP} or
a braided coaddition as in  \cite{M}) do not allow to get
braided vector fields on a q-hyperboloid. Also, note that our
braided vector fields differ drastically from  q-analogs of differential operators
arising from a coalgebraic structure in
the corresponding QG (cf. \cite{D}).

As for the q-Minkowski space algebra $\K_q[\R^4]$, it has one
generator more compared to $\K_q[\R^3]$. Moreover, this generator
is central. So, first, we define the partial derivative w.r.t.
this generator, and then we introduce the Maxwell operators on the
algebra $\K_q[\R^4]$ in the classical manner.

By properly defining the action of the QG $\uq$ on all ingredients
of the Maxwell operators on the algebras in question, we force
them to be
 $\uq$-invariant. Moreover, these operators  possess gauge freedom
 similar to the classical one (i.e. their kernels are as large as the kernels of their classical
counterparts are), provided the corresponding Laplace
 operators are central.

In conclusion we want to mention the  paper \cite{S} where the
author suggests a way of defining  q-analogs of gauge models
via quantum gauge potentials. For this end he
uses a q-analog of the Lie algebra $su(n)$ similar to that
considered above.
However, the ground (source) algebra considered in \cite{S} is commutative
whereas our ground algebras are essentially noncommutative.

We hope our method will be useful for a "q-deformation" of other gauge models.

\noindent{\bf Acknowledgement.} One of the authors (D.G.) would
like to thank the Max-Planck-Institut f{\"u}r Mathematik, where
this work was completed, for the warm hospitality and stimulating
atmosphere. The work is partially supported by the grant
ANR-05-BLAN-0029-01.

\section{Maxwell operator via projective modules}

In this section we introduce the Maxwell operator in classical settings and consider its behavior
w.r.t. the restriction to a subvariety.
Also, we give a few basic examples.

Let $\M$ be a regular affine variety endowed with a (pseudo-) Riemannian metric $g_{ij}=g(\partial_i,\partial_j)$
where $\partial_i$ are
partial derivatives in  local coordinates. We need two operators
$$ d: {\K[\M]} \to
{\Omega^1(\M)},\quad f(x)    \rightarrowtail \partial_i f\, dx_i,$$
$$\partial : {\Omega^1(\M)\to \K[\M]},\quad \al_i dx_i  \rightarrowtail {1\over \sqrt{g}} \partial_i(g^{ij}\sqrt{g}\, \al_j)$$
where $g=|\det(g_{ij})|$ and the tensor $g^{ij}$ is inverse to
$g_{ij}$. The Laplace  operator on  the algebra $\K[\M]$ is
$$\De(f)=\partial\, d f={1\over \sqrt{g}} \partial_i(g^{ij}\sqrt{g}\,\partial_j f).$$

 Laplace operators on the spaces $\Om^i(\M)$  are defined by
the formula
$$\De =\partial \,d+d\, \partial$$
where $\partial:{\Om^i(\M)}\to \Om^{i-1}(\M)$ is the well known
analog of the above operator. In what follows we realize the
Maxwell operator as $\Mw=\De-d\,\partial$. Besides, if ${\M}=\R^n$
and  the metric is constant in the Cartesian coordinates $x_i$,
the operator $\De$ acts on the space $\Om^1(\M)$ via $\De(\al_i
dx_i)= \De(\al_i) dx_i$.

\begin{proposition} Let $\N\subset \M$ be a regular subvariety of a variety $\M$ of codimension 1. Suppose that
in a vicinity of each point $a\in \N$ there exists a coordinate
system $(x_1,x_2,...,x_n)$ such that $\N$ is given by $x_n=0$ and
$(x_1,x_2,...,x_{n-1})$ is a local coordinate system in $\N$,
$g(\partial_n, \partial_n)=1$ and $g(\partial_i,
\partial_n)=0$ for any $1\leq i \leq n-1$. Then
$\Mw_{\N}=\Mw_{\M}|_{\N}$, i.e. the Maxwell operator on $\N$ is
the restriction of the Maxwell operator on $\M$. A similar claim
is valid for the Laplace operator.
\end{proposition}

{\bf Proof} follows from explicit form of the Maxwell and Laplace operators in the local
coordinate system $(x_1,x_2,...,x_n)$.

This proposition enables us to write  the Maxwell operator on
certain algebraic varieties in terms of ambient spaces. Thus, we
realize the space of one-forms (as well as the space of vector
fields) on such a variety as a projective module  without using
any local coordinate system.

Let us consider the basic example--a sphere
$$\sr=\{(x,y,z)\in \R^3\,|\, x^2+y^2+z^2=r^2,\, r>0\}$$
 embedded in the Euclidean space $\R^3\cong so(3)^*$ as an orbit of action of the group
$SO(3)$. Also,  we assume this space to be equipped with an
$SO(3)$-invariant pairing $\langle x, x \rangle=1,\,\, \langle x, y
\rangle=0,$ and so on. The corresponding metric is
$g(\partial_x,\partial_x)=1,\,\,g(\partial_x,\partial_y)=0,$ and
so on. Also, we endow the coordinate algebra  $\K[\R^3]$ of the
space $\R^3$ with $SO(3)$-covariant Poisson bracket
$$\{x,y\}=z, \,\,\{y,z\}=x, \,\,\{z,x\}=y.$$
To any function $f\in\K[\R^3] $ we associate the operator
$$\Pois_f(g):=\{f,g\},\,\, \forall\,g\in \K[\R^3].$$
Then the operators $X=\Pois_x,\,Y=\Pois_y,\, Z=\Pois_z$ are
infinitesimal rotations. Their explicit forms are
$$X=z\, \partial_y-y \, \partial_z,\,\,Y=x\, \partial_z-z\, \partial_x,\,\,Z=y\, \partial_x-x\, \partial_y.$$

They are tangent to the spheres $S^2_r$ and subject to the relation
$$x\,X+y\,Y+z\,Z=0.$$
Consider the coordinate algebra of the sphere $\sr$
$$\K[S^2_r]=\K[\R^3]/\langle x^2+y^2+z^2-r^2\rangle.$$
Hereafter $\langle I \rangle$ stands for the two-sided ideal
generated by a set $I$. The space $\Vect(S^2_r)$ of all vector
fields on a sphere (with coefficients from $\K[S_r^2]$), treated
as a $\K[S_r^2]$-module, is  the quotient
$$M=\K[S_r^2]^{\oplus 3}/{\overline M}$$
of the free $\K[S_r^2]$-module $\K[S_r^2]^{\oplus 3}$ over the submodule ${\overline M}=\{\varphi(x\,X+y\,Y+z\,Z),\,
\forall \,\varphi\in \K[S_r^2] \}$.

It is not difficult to see that the  module ${\overline M}$ is projective. Indeed, the matrix
$${\overline e}={1\over r^2}\left(\begin{array}{c}
x\\
y\\
z
\end{array}\right)
\left(\begin{array}{ccc}
x&y&z
\end{array}\right)$$
defines an idempotent such that ${\overline M}={\overline e}\,\K[S_r^2]^{\oplus 3}$. Therefore the
$\K[S_r^2]$-module $M$ can be realized as a submodule
$$M=e\,\K[S_r^2]^{\oplus 3}\subset \K[S_r^2]^{\oplus 3}$$
 where $e=I-{\overline e}$ is the complementary idempotent.

We call the $\K[S_r^2$]-module $M$ {\em tangent}.
 In contrast with other $\K[S_r^2]$-modules, the tangent module defines
the action $M\ot \K[S_r^2]\to \K[S_r^2]$  consisting in applying a
vector field to a function.

In a similar way we can realize the  space of one-forms
$\Om^1(S_r^2)$. This space consists of all elements $\{\al\,
dx+\beta\, dy+\gamma\, dz\}$ whereas the elements  ${\overline
M}=\{\varphi(x\,dx+y\,dy+z\,dz),\, \forall \,\varphi\in \K[S_r^2]
\}$ vanish. Thus, as a $\K[S_r^2]$-module, the space
$\Om^1(S_r^2)$ can be treated as the former module. The latter
module is called {\em cotangent}. (For a similar treatment of the
space $\Om^2[S_r^2]$  the reader is referred to \cite{GS1}.)

We do not distinguish between the tangent and cotangent
$\K[S_r^2]$-modules, and  denote them $M(S^2_r)$. Its elements are
treated to be triples $(\al,\,\beta,\,\gamma)^t$ ($t$ stands for
transposing) modulo
$$(x\,\rho,\, y\,\rho,\,z\,\rho)^t,\,\,\,\al,\,\beta,\,\gamma,\,\rho\in\K[S^2_r].$$

Now, we exhibit  the Maxwell operator on the Euclidian space $\R^3$ in a convenient form.
This operator acts on the space   of one-differential forms
$$\Om^1(\R^3)=\{\al\,dx+\beta dy+\gamma dz\},
\qquad \al,\beta,\gamma\in \K[\R^3]$$
(which is a free $\K[\R^3]$-module $\Om^1(\R^3)\cong \K[\R^3]^{\oplus 3}$) via formula (\ref{Max}).

It can be also written in the following form
$$-\rot\,\rot=\De-\grad\, \div$$
where
$$\rot:\Om^1(\R^3)\to \Om^1(\R^3), \, \div:\Om^1(\R^3)\to \K[\R^3],
\,\grad:\K[\R^3]\to \Om^1(\R^3)$$
are the curl, divergence, and gradient respectively and
$$\De=\De_{\K[so(3)^*]}=\partial^2_x+\partial^2_y+\partial^2_z$$
is  the Laplace operator.

By  identifying a differential
form $\al\,dx+\beta\,dy+\gamma\,dz$  and the triple $(\al,\,\beta,\,\gamma)^t$ as explained above  we can write
\be
\Mw_{\K[so(3)^*]}\left(\begin{array}{c}
\al\\
\beta\\
\gamma
\end{array}\right)=
\left(\begin{array}{c}
\De(\al)\\
\De(\beta)\\
\De(\gamma)
\end{array}\right)-\left(\begin{array}{c}
\partial_x\\
\partial_y\\
\partial_z
\end{array}\right)(\partial_x,\,\partial_y,\,\partial_z)\left(\begin{array}{c}
\al\\
\beta\\
\gamma
\end{array}\right).
\label{Mw} \ee

 Observe that the Maxwell operator is $SO(3)$-invariant  if $SO(3)$ acts on the generators $x,y,z$  and the matrices
 $M\in {\Mat}_3(\R)$ in a proper way (see section 5).
The gauge freedom is due to the fact that the triples $(\partial_x\rho, \,
\partial_y\rho,\, \partial_z\rho)^t$ belong to the kernel $\Ker\,(Mw_{\K[so(3)^*]})$ of this operator. So, if a triple
$(\al,\,\beta,\,\gamma)^t$ is a solution of the Maxwell equation
$\Mw_{\K[so(3)^*]}(\al,\,\beta,\,\gamma)^t=(\lambda,\,\mu,\,\nu)^t$
then the triple $(\al,\,\beta,\,\gamma)^t+(\partial_x\rho, \,
\partial_y\rho,\, \partial_z\rho)^t$ is also a solution.

Now, consider the Maxwell operator on the sphere $\sr$. By using
the relation \be r\partial_r=x\partial_x+y\partial_y+z\partial_z
\label{dirr} \ee we realize
 the partial derivatives $\partial_x,\,\partial_y,\,\partial_z$ as follows
\be
\partial_x={1\over r^2}(y\, Z-z\, Y)+{x\over r}\partial_r \, \,\circlearrowleft. \label{part}
\ee
(Thereafter the symbol $\circlearrowleft$ stands for the cyclic permutations).
In what follows we also use the following formulae
$$\partial_r x={x\over r} \,\, \circlearrowleft,\,\,\partial_x r={x\over r}\,\, \circlearrowleft,\,\,\,X(r)=Y(r)=Z(r)=0,$$
and the fact that the vector field $\partial_r$ commutes with
$X,Y,Z$.

By using formula (\ref{part}) we rewrite the Laplace operator on the space $\R^3$ as follows
\be
\De_{\K[so(3)^*]}=({1\over r^2}\X+{x\over r}\partial_r)^2+ ({1\over r^2}\Y+{y\over r}\partial_r)^2+
({1\over r^2}\Z+{z\over r}\partial_r)^2, \label{Lap}\ee
where we use the notations
$$\X=y\, Z-z\, Y,\qquad \Y=z\, X-x\, Z,\qquad \Z=x\, Y-y\, X.$$

By  the above proposition the Laplacian $\De_{\K[S^2_r]}$ on the
sphere $S^2_r$
 is  the restriction of the Laplacian $\De_{\K[so(3)^*]}$ to the sphere $S^2_r$. Indeed, in  spherical coordinates
the radius $r$ plays the role of the coordinate $x_n$ from Proposition 1.

\begin{lemma}
\be
\De_{\K[S^2_r]}={\X^2+\Y^2+\Z^2\over r^4}. \label{Lapl} \ee
\end{lemma}

{\bf Proof} We have only to check that the first
order component of the operator (\ref{Lap}) vanishes on the sphere
 $S^2_r$ $\partial_r=0$. Indeed, this component equals
$${x\over r}\partial_r({1\over r^2})\X+{y\over r}\partial_r({1\over r^2})\Y+{z\over r}\partial_r({1\over r^2})\Z
=-{2\over r^4}(x\X+y\Y+z\Z)=0.$$ (Moreover, there are no
$SO(3)$-invariant first order differential operators on the space
$\R^3$.)
\begin{lemma}
The Laplace operator (\ref{Lapl}) equals
$$\De_{\K[S^2_r]}={X^2+Y^2+Z^2\over r^2}.$$
\end{lemma}

{\bf Proof} We have
$$\X^2+\Y^2+\Z^2=(y\,Z-z\,Y)^2+(z\,X-x\,Z)^2+(x\,Y-y\,X)^2=$$ $$y^2Z^2+z^2Y^2-yz(YZ+ZY)-yxZ+zxY+\circlearrowleft$$ $$
=(r^2-x^2)X^2-yxYX-zxZX+\circlearrowleft=$$ $$r^2(X^2+Y^2+Z^2)-[x(xX+yY+zZ)X+\circlearrowleft]=r^2(X^2+Y^2+Z^2).$$

This form of the Laplacian is more familiar and widely used in study of  rotationally  symmetric
Schrodinger operators.

Let us emphasize that the operators $X,\,Y,\,Z,\,\,\X,\, \Y,\,Z$
are well defined on the space $\R^3$, and the relation
$\X^2+\Y^2+\Z^2=r^2(X^2+Y^2+Z^2)$ is also valid on the whole space
$\R^3$.

\begin{proposition} \label{th1}
$$ \De_{\K[\sr]}\,\X-\X \,\De_{\K[\sr]}=-{{2\,x}}\, \De_{\K[\sr]}\,\,\circlearrowleft$$
on $\R^3$.
\end{proposition}
{\bf Proof} Indeed, by direct computations we have
$$\De_{\K[S^2_r]}\,\X-\X\, \De_{\K[S^2_r]}={2\over r^2}(zZX+yYX-xY^2-xZ^2)={2\over r^2}((xX+yY+zZ)X-x(X^2+Y^2+Z^2)).$$

Now,  define the Maxwell operator  $\Mw_{\K[S^2_r]}$  on the
sphere $S_r^2$ as follows
 \be
 \Mw_{\K[S^2_r]}\left(\begin{array}{c}
\al\\
\beta\\
\gamma
\end{array}\right)=e\left(
\left(\begin{array}{c}
\De_{\K[S^2_r]}(\al)\\
\De_{\K[S^2_r]}(\beta)\\
\De_{\K[S^2_r]}(\gamma)
\end{array}\right)-{1\over r^4}\left(\begin{array}{c}
\X\\
\Y\\
\Z
\end{array}\right)
(\X,\,\Y,\,\Z)\left(\begin{array}{c}
\al\\
\beta\\
\gamma
\end{array}\right)\right),\,\,\,\,\,\,
\left(\begin{array}{c}
\al\\
\beta\\
\gamma
\end{array}\right)\in M(S^2_r).
\label{op} \ee

In this definition we assume that the elements of the module
$M(\sr)$ are triples $(\al,\,\beta,\,\gamma)^t$ such that
${\overline e}\,(\al,\,\beta,\,\gamma)^t=0$ or, which is the same,
$ e\, (\al,\, \beta,\,\gamma)^t=(\al,\,\beta,\,\gamma)^t$. The
idempotent  $e$ coming in this definition ensures that the image
of the operator $\Mw_{\K[S^2_r]}$ belongs to the module $M(\sr)$.

Now, we justify this definition. In fact, by  Proposition 1 if we
present the Maxwell operator on $\R^3$ in spherical coordinates
and restrict it to the sphere $\sr,$ we get the Maxwell operator
on this sphere. It remains to check that it coincides with the
operator (\ref{op}).

Let us observe that similarly to the Laplacian (\ref{Lapl}) the
operator  $\Mw_{\K[S^2_r]}$ is $SO(3)$-covariant. Moreover, its
gauge freedom consists in the fact that the triples
$(\X(\rho),\,\Y(\rho),\,\Z(\rho))^t$ belong to
$\Ker\,(\Mw_{\K[S^2_r]})$. It follows from the proposition
\ref{th1}. (These triples belong to the module $M(\sr)$ since
${\overline e}(\X,\,\Y,\,\Z)^t=0$.)

\begin{remark} Introduce "gradient" and "divergence" on the sphere $\sr$ as follows
$$\grad_{\K[S^2_r]}\, f=r^{-2}(\X(f),\,\Y(f),\,\Z(f)),\qquad \div_{\K[S^2_r]}\left(\begin{array}{c}
\al\\
\beta\\
\gamma
\end{array}\right)=r^{-2}(\X,\,\Y,\,\Z)\left(\begin{array}{c}
\al\\
\beta\\
\gamma
\end{array}\right),\,\,\,\,\,\,
\left(\begin{array}{c}
\al\\
\beta\\
\gamma
\end{array}\right)\in M(S^2_r).$$
Rewrite the operator $\Mw_{\K[S^2_r]}$ in the form similar to the
classical one:
$$\Mw_{\K[S^2_r]}=e\,\De_{\K[S^2_r]}-\grad_{\K[S^2_r]}\,\div_{\K[S^2_r]}$$
(the factor $e$ in the second summand can be omitted).
\end{remark}

Let us consider one example more--the classical Minkowski space, i.e. 4 dimensional space endowed  with an
$SO(1,3)$-covariant  norm
$$ \|(t,\,x,\,y,\,z )\|= \sqrt{t^2-x^2-y^2-z^2}.$$
The corresponding second order differential operator is
$$\De_{\K[\R^4]}=\partial^2_t-\partial^2_x-\partial^2_y-\partial^2_z.$$

It is called d'Alembertian or $so(1,3)$ Laplacian.
Then the corresponding Maxwell operator is
\be
\Mw_{\K[\R^4]}\left(\begin{array}{c}
\al\\
\beta\\
\gamma\\
\delta
\end{array}\right)=
\left(\begin{array}{c}
\De_{\K[\R^4]}(\al)\\
\De_{\K[\R^4]}(\beta)\\
\De_{\K[\R^4]}(\gamma)\\
\De_{\K[\R^4]}(\delta)
\end{array}\right)-\left(\begin{array}{c}
\partial_t\\
\partial_x\\
\partial_y\\
\partial_z
\end{array}\right)(\partial_t\,-\partial_x,\,-\partial_y,\,-\partial_z)\left(\begin{array}{c}
\al\\
\beta\\
\gamma\\
\delta
\end{array}\right).
 \ee

It is evident that it is $SO(1,3)$-covariant and
$$(\partial_t\rho,\, \partial_x\rho,\,\partial_y\rho,\,\partial_z\rho)^t\in \Ker\, (\Mw).$$

It is also clear that by restricting $t=0$ we get the Maxwell
operator on the  space $\R^3$ (up to a sign).

\section{Maxwell operator on $sl(2)^*$ and $\Hr$}

Now, apply the above scheme to the space $\R^3\cong sl(2)^*$
endowed with an action of the group $SL(2)$. This example is
another (non-compact) real form of the situation considered above.
So, the corresponding Maxwell operator can be obtained by a mere
change of a basis. Nevertheless, we describe it in detail since it
is going to  be "q-deformed" below.

Consider a basis $\{b,\, h,\,c\}$ in the space $\K[sl(2)^*]$
equipped with the $SL(2)$-covariant Poisson bracket
$$\{h,b\}=2b,\quad \{h,c\}=-2c,\quad \{b,c\}=h.$$
The corresponding Poisson operators are
$$H={\Pois}_h=2b\,\partial_b-2c\,\partial_c,\quad B={\Pois}_b=h\,\partial_c-2b\, \partial_h, \quad C={\Pois}_c=-h\,\partial_b+2c\,\partial_h.$$
They are tangent to hyperboloids and subject to the relation
\be
c\,B+{hH\over 2}+b\,C=0. \label{tang} \ee

 The $\K[\Hr]$-module $\Vect(\Hr)$
of vector fields on a hyperboloid (\ref{hyper}) is a quotient
module $M=\K[\Hr]^{\oplus 3}/{\overline M}$ of the free
$\K[\Hr]$-module $\K[\Hr]^{\oplus 3}$ over the submodule
${\overline M}=\{\varphi(c\,B+{hH\over 2}+b\,C),\, \forall
\,\varphi\in \K[\Hr] \}$. The idempotent corresponding to the
module ${\overline M}$ is \be {\overline e}={1\over
r^2}\left(\begin{array}{c}
c\\
{h\over 2}\\
b
\end{array}\right)
\left(\begin{array}{ccc}
b&h&c
\end{array}\right)
\label{idem} \ee (In order to show that ${\overline M}={\overline
e}\K^{\oplus 3}$  it  suffices to check that
$\K[\Hr]=\{b\al+h\beta+c\gamma\}$.) Thus, the module $M$ is also
projective $M=e\,\K[\Hr]^{\oplus 3}$ where $e=1-{\overline e}$.
The $\K[\Hr]$-module $\Om^1(\Hr)$ can be treated similarly.

 Endow the space ${\span}(b,h,c)$ with an $SL(2)$-covariant pairing
 \be
 \langle b,\,c \rangle=\langle c,\,b \rangle=1,  \langle h,\,h \rangle=2 \ee
 which is inverse to the  Casimir element
 $${\Cas}=bc+{h^2\over 2}+cb.$$
  Thus, on the space $sl(2)$
 \be
 \langle b=\partial_c,\quad \langle c=\partial_b,\quad \langle h=2\partial_h \label{equ} \ee
where $\langle x:sl(2)\to \K$ is the "bra" operator, i.e., such that $\langle x\,(y):=\langle x,\, y\rangle$.

We extend the action of the operators $\langle b,\, \langle c,\,
\langle h$ to the algebra $\K[sl(2)^*]$ via the relations
(\ref{equ}), i.e. by means of the Leibnitz rule. Thus, the action
$$sl(2)\ot \K[sl(2)^*]\to \K[sl(2)^*]$$
is well defined. It is clear that it is $SL(2)$-covariant.
Otherwise stated, we have an $SL(2)$-covariant map \be
sl(2)\to\Vect(sl(2)^*)\label{map}\ee different from that defined
above via the Poisson bracket. We have associated a tangent
(Poisson) vector field to any element from $sl(2)$, now it is a
partial derivative which is associated to such an element.

By using the map (\ref{map}) we associate a differential operator
to any element from $U(sl(2))$. Thus, the Casimir element is
related to the following operator
$$\De_{\K[sl(2)^*]}=\partial_b\partial_c+{2\partial_h^2}+\partial_c \partial_b.$$
It is a non-compact analog of the  Laplace operator $\De_{\K[so(3)^*]}$.  (More precisely, it is a multiple of the latter
Laplacian  written in the basis $\{b,\, h,\,c\}$.)

Similarly, the element $\Id\,\Cas-{\overline e}$ is related to the
operator \be \Mw_{\K[sl(2)^*]}\left(\begin{array}{c}
\al\\
\beta\\
\gamma
\end{array}\right)=
\left(\begin{array}{c}
\De(\al)\\
\De(\beta)\\
\De(\gamma)
\end{array}\right)-\left(\begin{array}{c}
\partial_b\\
\partial_h\\
\partial_c
\end{array}\right)(\partial_c,\,{2}\partial_h,\,\partial_b)\left(\begin{array}{c}
\al\\
\beta\\
\gamma
\end{array}\right),\,\,\left(\begin{array}{c}
\al\\
\beta\\
\gamma
\end{array}\right)\in \K[sl(2)^*]^{\oplus 3}.
 \ee
Note, that  $(\partial_b \rho,\, \partial_h \rho,\,\partial_c \rho)^t \in \Ker(\Mw_{\K[sl(2)^*]})$.

Now, we want to relate the vector fields $\langle b,\,\langle
h,\,\langle c$ to the fields tangent to all hyperboloids. Observe
that the following formula, similar to (\ref{dirr}), holds:
$$r\, \partial_r=b\,\partial_b+h\,\partial_h+c\,\partial_c.$$
This entails the following relations, similar to (\ref{part}) \be
\langle b=\partial_c={h\, B-b\, H \over 2 r^2}+{b \over
r}\partial_r, \quad {1\over 2}\langle h=\partial_h={b\, C-c\, B
\over 2 r^2}+{h \over 2r}\partial_r, \quad \langle
c=\partial_b={c\, H-h\, C \over 2 r^2}+{c \over r}\partial_r .
\label{rela} \ee Introduce the following notations
$$\B={1\over 2}(h\, B-b\,H),\quad \H=b\, C-c\, B,\quad \C={1\over 2}(c\, H-h\, C).$$
Thus, we get
\be
\langle b={\B\over r^2}+{b\over r}\partial_r,\quad \langle h={\H\over r^2}+{h\over r}\partial_r,\quad
\langle c={\C\over r^2}+{c\over r}\partial_r. \label{raven}
\ee

Now,  introduce Laplace operator on a hyperboloid $\Hr$ similarly
to the previously described:
$$\De_{\K[\Hr]}={1\over r^4}(\B\, \C+{\H^2\over 2}+\C\, \B).$$
Also,
$$\De_{\K[\Hr]}={1\over r^2}(B\, C+{H^2\over 2}+C\, B). $$

According the  pattern above we define the Maxwell operator on the
hyperboloid $\Hr$ as
$$\Mw_{\K[\Hr]}\left(\begin{array}{c}
\al\\
\beta\\
\gamma
\end{array}\right)=e\left(
\left(\begin{array}{c}
\De_{\K[\Hr]}(\al)\\
\De_{\K[\Hr]}(\beta)\\
\De_{\K[\Hr]}(\gamma)
\end{array}\right)-{1\over r^4}\left(\begin{array}{c}
\C\\
{\H\over 2}\\
\B
\end{array}\right)
(\B,\,\H,\,\C)\left(\begin{array}{c}
\al\\
\beta\\
\gamma
\end{array}\right)\right),\,\,\,\,\,\,
\left(\begin{array}{c}
\al\\
\beta\\
\gamma
\end{array}\right)\in M(\Hr).
$$

It is evident that the triples $(\C\rho,\, {\H\rho \over 2},\, \B\rho)^t$ belong
to the kernel  $\Ker\,(\Mw_{\K[\Hr]})$.

Completing this section we want to emphasize that the Maxwell
operators on the algebras $\K[\R^3]$ and $\K[\Hr]$ are
$SL(2)$-invariant provided the matrices coming in the definition of these
operators are endowed with a proper action of the group $SL(2)$.
We exhibit a way to introduce such an action in the last section
in a more general settings of quantum algebras.

\section{Elements of analysis   on q-Minkowski space algebra}

As it has been mentioned above,  the q-Minkowski space algebra is
a particular case of the Reflection Equation Algebra, which is
defined as follows. Let $R:\vv\to\vv $ be a {\em braiding}, i.e.,
an invertible operator subject to the quantum Yang-Baxter equation
$$(R\ot I)(I\ot R)(R\ot I)=(I\ot R)(R\ot I)(I\ot R)$$
where $V$ is a $n$-dimensional ($n\geq 2$) vector space over the field $\K$.

Let $L=\|l_i^j\|$ be a $n\times n$ matrix with entries
$l_i^j,\,\,1\leq i, j \leq n$. Then the relation \be R(L\ot
I)R(L\ot I)-(L\ot 1)R(L\ot I)R=0 \label{RE} \ee is called the {\em
reflection equation}. The algebra generated by the unity and the
elements $l_i^j$ subject to this system is called a {\em
reflection equation algebra} (REA) and denoted $\lrq$.

If, in addition,  the braiding $R$  is  subject to the relation
$$(qI-R)(\qq I+R)=0, \,\,q \in \K$$
it is called a {\em Hecke symmetry}. If $q=1$ it becomes an {\em
involutive symmetry}.

Let $n=2$ and  $R$ be the product of the image in the space $\vv$
of the universal R-matrix of the QG $\uq$ and the usual flip. Then
in the basis $\{x_1\ot x_1, x_1\ot x_2, x_2\ot x_1, x_2\ot x_2\}$
(where $\{x_1,\, x_2\}$ in an appropriate basis in $V$) the
braiding $R$ reads \be R_q=\left(\begin{array}{cccc}
q&0&0&0\\
0&q-\qq&1&0\\
0&1&0&0\\
0&0&0&q
\end{array}\right).
\label{canon}
\ee
It is easy to see that $R$ is a Hecke symmetry. The parameter $q$ is assumed to be generic.

In what follows the algebra $\lrq$ corresponding to this Hecke symmetry is
called q-{\em Minkowski space algebra}. In this case
it is also denoted $\K_q[\R^4]$.
Also, we use the following notation for generators of this algebra
$$L=\left(\begin{array}{cc}
l_1^1&l_1^2\\
l_2^1&l_2^2 \end{array}\right)=
\left(\begin{array}{cc}
a&b\\
c&d \end{array}\right).$$

Let us explicitly write down the system (\ref{RE}) with the Hecke
symmetry (\ref{canon}): \be
\begin{array}{l@{\hspace{20mm}}l}
q ab - q^{-1}ba = 0 &q( bc - cb)  =(q-\qq) a(d-a)\\
q ca - q^{-1}ac  = 0 & q(cd - dc)  = (q-\qq)c a\\
ad-da=0 &  q(db - bd)   = (q-\qq) a b.
\end{array}
\label{sys1}
\ee

Now, rewrite  the system (\ref{sys1}) in the basis
$\{l,\,h,\,b,\,c\}$ where $l=\qq a+q d,\quad h=a-d$: \be
\begin{array}{l@{\hspace{20mm}}l}
q^{2}hb-bh = -(q-\qq) lb & bl = lb\\
q^{2}ch-hc  =-(q-\qq) lc & cl = lc\\
\left(q^{2}+1\right)\left(bc-cb\right)+\left(q^{2}-1\right)h^{2} = -(q-\qq) lh& hl=lh
\end{array}
\label{sys2}
\ee

Observe that the element $l$ is central but, opposite to the
classical case, it appears in the equations of the left column of
the system (\ref{sys2}). Also, we need the algebra
$\K_q[\R^3]=\K_q[\R^4]/\langle l \rangle$ which is a braided
analog of the coordinate algebra $\K[\R^3]$. It is generated by
three elements $b,\, h,\, c$ subject to \be q^{2}hb-bh
=0,\,\,\,\,q^{2}ch-hc  =0,\,\,\,\, (q^{2}+1)(bc-cb)+(q^{2}-1)h^{2}
= 0. \label{defi} \ee

The generating spaces of this algebras ${\span}(a,b,c,d)={\span}(b,h,c,l)$ and ${\span}(b,h, c)$
are respectively denoted $\LL$ and $\SL$.
Let us endow them with an action of the QG $\uq$ as follows.

Recall that the QG $\uq$ is generated by the unit and four
generators $K,\,K^{-1},\,X,\,Y$ subject to the relations
$$K\,K^{-1}=1,\quad K^\epsilon\, X=q^{2\epsilon} X\, K^\epsilon\, \quad
K^\epsilon\,Y=q^{-2\epsilon} Y\, K^\epsilon,\quad XY-YX={{K-K^{-1}}\over{q-\qq}},\quad \epsilon\in\{-1,1\}.$$

There exists a family of coproducts and corresponding antipodes which (together with the standard counit)
 endow this algebra with a Hopf structure. This family is parameterized by a continuous parameter $\theta$
 assumed to be a real number (cf. \cite{DG}):
$$\De(K^\epsilon)=K^\epsilon \ot K^\epsilon,\quad \De(X)=X\ot K^{\theta-1}+K^\theta\ot X,\quad
\De(Y)=Y\ot K^{-\theta}+K^{1-\theta}\ot Y.$$
(In fact, all coproducts are equivalent, cf. \cite{DG}.)

We define an action of the QG $\uq$ on the space $\SL$ as follows
$$K^\epsilon(b)=q^{2\epsilon} b,  \,\,K^\epsilon(h)=h,  \,\,     K^\epsilon(c)=q^{-2\epsilon}c,$$
$$X(b)=0,  \,\,X(h)=-q^{\theta}2_q b,\,    X(c)= q^{1-\theta}h,$$
$$Y(b)=-q^{-\theta}h,  \,\,Y(h)=q^{\theta-1}2_q c,\,\,    Y(c)=0.$$
Hereafter, $n_q={q^n-q^{-n}\over q-\qq}$. Moreover, we put
$X(l)=Y(l)=0,\, K^\epsilon(l)=l$. Thus, as a $\uq$ module, the
space $\LL$ is a direct sum of  $\SL$ and $\K\, l$.

The reader can easily check that the structure of the  algebra
$\K_q[\R^3]$ is compatible with the action of the QG $\uq$
extended to higher powers of $\SL$ via the coproduct $\De$. In
order to do this it suffices to check that the system (\ref{defi})
is invariant w.r.t. the QG $\uq$.

Now, describe a regular way of introducing Casimir-like elements
and covariant pairing into these spaces.
 In order to do this we need two special operators playing the role
of the parity operators in super-spaces.

Let $R:\vv\to\vv$ be a braiding. It is  called {\em skew-invertible} if there exists an operator
$\Psi:\vv\to\vv$ such that
$${\rm Tr}_{2}R_{12}\Psi_{23}=P_{13}$$
where $P_{13}$ is the flip transposing  the first and third spaces in the product $V^{\ot 3}$
and ${\rm Tr}_{2}$ stands for the trace applied in the second  space.

Fixing a basis $\{x_i\}\in V$  and representing the operators $R$
and ${\Psi}$  in the basis  $\{x_i\ot x_j\}\in \vv$
  by the matrices $R_{ia}^{kb}$ and  ${\Psi}^{al}_{bj},$ respectively, we  rewrite this relation
as follows
$$
\sum_{a,b = 1}^n R_{ia}^{kb} {\Psi}^{al}_{bj} = \delta^l_i\,
\delta^k_j.
$$

For any skew-invertible Hecke symmetry the operators
$$B:={\rm Tr}_{1}\Psi\,\,(B_i^j={\Psi}^{aj}_{ai}),\,\,\,\,C:={\rm Tr}_{2}\Psi\,\,(C_i^j={\Psi}^{ja}_{ia})$$
are well defined and  the elements $\Tr_q \, L^k=\Tr\, C\, L^k$ are central in the algebra $\lrq$ (cf. \cite{GPS1}).

The operator  $L^k\to \Tr_q \, L^k$ is called the {\em quantum
trace}. It can be  extended by linearity to all polynomials of the
matrix $L$. For the q-Minkowski space algebra $\K_q [\R^4]$ we
have $C={\diag}(q^{-3},\,\qq)$. Hence, $l=q^2\,\Tr_q \, L$. Also,
 we consider the central element $q^2\,\Tr_q\,
L^2$ in this algebra. In the basis $\{b,c,h,l\}$ its explicit form
is \be {\Cas}_{gl}=\qq bc +q cb+{1\over{2_q}} (h^2+l^2).
\label{casgl} \ee

Its image in the algebra $\K_q[\R^3]$ is \be {\Cas}_{sl}=\qq bc
+{1\over{2_q}} h^2+q cb.  \label{cas} \ee It is a central element
in this algebra. So, we get braided analogs of the $gl(2)$ and
$sl(2)$ Casimir elements, respectively. The quotient of the
algebra $\K_q[\R^3]$ over the ideal generated by the element
${\Cas}_{sl}-r^2$ is called a {\em quantum (braided or
q-)hyperboloid}.

\begin{remark} Observe that ${\Cas}_{sl}$ is the unique quadratic element
(up to a factor) in $\K_q[\R^3]$ which is $\uq$-invariant. As for the algebra
$\K_q[\R^4]$, since  the element $l^2$ is also $\uq$-invariant we have a
family of such elements, namely, all linear combinations
of ${\Cas}_{sl}$ and $l^2$.
\end{remark}

We also need $\uq$-invariant pairings in the spaces $\LL={\span}(b,c,h,l)$ and $\SL={\span}(b,c,h)$.
For a general skew-invertible Hecke symmetry
such a pairing on the space $\LL={\span}(l_i^j)$ is defined via the operator $B$. Namely, we put
$$\langle l_i^j,\, l_k^l\rangle=B_k^j \delta_i^l.$$
Also, an analog of the space $\SL$ can be defined in this case as
the subspace of $\LL$ which is the kernel of the linear map
defined on the generators $l_i^j$ by $l_i^j \to \delta_i^j$. Note,
that the space $\LL$ can be identified with  $\End(V)$ so that
this map is nothing but a braided analog of the trace written in
the basis $\{l_i^j\}$  (cf. \cite{GPS2}).

In the case  of the q-Minkowski space algebra the pairing table is
(we only exhibit terms with non-trivial pairing): \be \langle
a,a\rangle=\qq,\,\, \langle b,c \rangle=q^{-3},\,\, \langle c,b
\rangle=\qq,\,\, \langle d, d \rangle=q^{-3}\,\,{\rm hence}\,\,
\langle h, h\rangle=q^{-2}2_q,\,\, \langle l,l
\rangle=q^{-2}2_q.\label{pair} \ee So, contrary to the classical
case, this pairing is not symmetric. Note, that the restriction of
this paring to the space $\SL$ is the unique (up to a factor)
$\uq$-covariant pairing. However, we have a freedom for the
pairing $\langle l, l \rangle$ on the space $\LL$.

Our next goal is to define a braided analog of tangent vector fields. In order to do so, we endow
 the space $\SL$ with a braided analog of the Lie bracket $sl(2)$. Such analogs
 of Lie brackets  $gl(n)$ and $sl(n)$ can be defined in the spaces $\LL$ and $\SL,$ respectively, for any skew-invertible
 Hecke symmetry $R$. We refer the reader to \cite{GPS2} for their construction (also, see remark \ref{rem}).
 But in the low dimensional case in question we define such a bracket
 on the space $\SL$   by only using the fact that it is covariant.

Let us extend the action of the QG $\uq$ to the space $\SL\ot\SL,$
and decompose it in a direct sum of irreducible $\uq$ submodules
$\SL\ot\SL=V_0\oplus V_1 \oplus V_2$ where the subscript stands
for the spin. Then the operator
$$[\,,\,]: \SL\ot \SL\to \SL$$
is a $\uq$ morphism iff it is  trivial on the components $V_0$ and
$V_2,$ and it is an isomorphism between $V_1$ and $\SL$. By this
condition the bracket is defined in the unique  (up to a factor)
way.

The multiplication table is as follows
$$ [b,b]=0,\quad [b,h]=-w \,b,\quad [b,c]=w \,{q \over 2_q}\,h,\quad [h,b]=w\,q^2\, b,$$
$$[h,h]=w\,(q^2-1),\quad [h,c]=-w\, c, \quad  [c,b]=-w\,{q\over 2_q}\, h,\quad [c,h]=w\, q^2\, c,\quad [c,c]=0.$$
Here $w$ is an arbitrary factor. It can be fixed if we introduce the "enveloping algebra"
of this "q-Lie algebra"  by the relations
\be q^{2}hb-bh =\h b,\,\,\,\,(q^{2}+1)(bc-cb)+(q^{2}-1)h^{2} = \h h,\,\,\,\, q^{2}ch-hc  =\h c   \label{env} \ee
and require the above bracket to define a representation of this enveloping algebra. Then $w=\h (q^4-q^2+1)^{-1}$.
We denote the space $\SL$ endowed with the above bracket $sl(2)_q$ and the algebra defined by
the relations (\ref{env}) $U(sl(2)_q)$.

Introduce q-analogs of the adjoint operators as follows
 $B_q=\ad(b),\,H_q=\ad(h),\, C_q=\ad(c)$ where the action $\ad$ is defined via the above bracket.
 These operators in the basis $\{b,h,c\}$ are
\be
B_q=w\left(\begin{array}{ccc}
0&-1&0\\
0&0&{q\over 2_q}\\
0&0&0
\end{array}\right)\quad
H_q=w\left(\begin{array}{ccc}
q^2&0&0\\
0&q^2-1&0\\
0&0&-1
\end{array}\right)
\quad
C_q=w\left(\begin{array}{ccc}
0&0&0\\
-{q\over 2_q}&0&0\\
0&q^2&0
\end{array}\right)\quad
\label{q-ad} \ee

\begin{proposition} The operators $B_q, \,H_q,\,C_q$ are subject to
\be
\qq b\,C_q+{h\,H_q\over 2_q}+q c\, B_q=0. {\rm \label{re}} \ee
\end{proposition}

{\bf Proof} is straightforward.

Now, we want to extend the operator $B_q,\,H_q,\,C_q$ to the
algebra $\K_q[\R^3],$  preserving relation (\ref{re}). It will be
done by a slightly modified  method suggested in \cite{A}.

Let $V_k$ be the $\uq$ submodule of $\SL^{\ot k}$ with the highest vector $b^k$, i.e.,
$$V_k={\span}(b^k,\,Y(b^k),\, Y^2(b^k),..., Y^{2k}(b^k)).$$
Note, that $\dim V_k=2k+1$.
There exists a $\uq$-covariant projector $P_k:\SL^{\ot k}\to V_k$. It can be realized as a polynomial in
$$\RR_{12}=\RR\ot  I_{2k-1},
\, \RR_{23}=I\ot \RR\ot I_{2k-2},...,\RR_{2k\,2k+1}=I_{2k-1}\ot
\RR$$ where $\RR$ is the product of  the  universal quantum
R-matrix  represented in the space $\SL^{\ot 2}$ and the flip. A
way of constructing such operators $P_k$ is described in
\cite{OP}. Observe, that the braiding $\RR$ is of the
Birman-Murakami-Wenzl type, and therefore the results of this
paper can be applied.  Then the extension of the operators
$B_q,\,H_q,\,\C_q$ to the component $V_k$ (denoted
$B_q^{(k)},\,H_q^{(k)},\,\C_q^{(k)},$
 respectively) are defined as follows
$$B_q^{(k)}=\tau_k\, P_k(B_q\ot I_{2k}),\,\,\,H_q^{(k)}=\tau_k\, P_k(H_q\ot I_{2k}),\,\,\,C_q^{(k)}=
\tau_k\, P_k(C_q\ot I_{2k})$$
where the factor $\tau_k$ can be found from the property that the
operators $B_q^{(k)},\,H_q^{(k)},\,\C_q^{(k)}$ realize a representation of the
the algebra $U(sl(2)_q)$. Thus, the prolongation of the operators $B_q,\,H_q,\,\C_q$ is
well defined on all components $V_k$.

Now, observe that for a generic $q$ we have $\K_q[\R^3] \cong(\oplus V_k)\ot Z$
where $Z$ is the center of the algebra $\K_q[\R^3]$.
By putting $B_q(v\ot z)=B_q^k(v)\ot z$ where $v\in V_k,\, z\in Z$ and analogously
for $H_q$ and $C_q$ we define the operators
$B_q, \,H_q,\,C_q$ on the whole algebra $\K_q[\R^3]$. We call them
{\em generating braided tangent vector fields}.
By considering their combinations with coefficients from
$\K_q[\R^3]$ we get, by definition, all {\em braided tangent vector fields}.
This definition is motivated by the following.

\begin{proposition} The extended operators $B_q, \,H_q,\,C_q$ are subject to the relation (\ref{re}).
\end{proposition}

Thus, the space of all braided tangent vector fields is treated as
a left $\K_q[\R^3]$ module which is the quotient of the free
module $\K_q[\R^3]^{\oplus 3}$ over the submodule generated by the
l.h.s. of (\ref{re}).

\begin{remark} \label{rem} Note that we have defined braided tangent vector fields without
using any form of the Leibnitz rule (in other words, any coproduct
in the algebra $\K_q[\R^3]$). Nevertheless, if $R$ is a
skew-invertible Hecke symmetry such a coproduct exists in the
so-called modified REA which is a braided analog of the enveloping
algebra $U(gl(n))$. This modified REA can be obtained from the REA
$\lrq$ by a shift (cf. \cite{GPS2} for details). It is possible to
use this coproduct in order to define braided analogs of tangent
vector fields in the space $\SL$ for $n>2$. However, it is not
clear whether vector fields satisfy the relations looking like
those of (\ref{re}). Nevertheless,
 we want to point out the property of the mentioned coproduct
to be a $\Uq$ morphism. By using the construction of the paper \cite{LS} it
is also possible to get a series of representations of the algebra
$U(sl(2)_q)$ close to those above.
But this way of construction the representations of this
 algebra does not give rise to braided tangent vector fields since the relation (\ref{re}) is not satisfied.
\end{remark}

\section{Braided Maxwell operator on quantum algebras}

In this section  we define the Maxwell operator on the algebras  $\K_q[\R^3]$, $\K_q[\Hr]$, and $\K_q[\R^4]$
 following the classical pattern discussed above.
First, we  introduce braided analogs of the relations
(\ref{rela}). In order to do this, we need to define braided
analogs of the operators $\B,\,\H,\,\C$. We put
$$\B_q=w^{-1} (q^2 hB_q-b H_q),\,\,\H_q=w^{-1} ((q^2+1)(bC_q-cB_q)+(q^2-1)hH_q) ,\,\,\C_q=w^{-1} (q^2cH_q-hC_q).$$
Here $w$ is a factor coming in (\ref{q-ad}). For the sake of
convenience  we put $w=2_q$. Then for $q\to 1$ we retrieve the
classical operators $\B,\,\H,$ and $\C$ respectively.

On the second step we define the derivative $\partial_r$ in
the algebra $\K_q[\R^3]$ in the classical way, i.e., we assume
  the derivation in $r$ to be subject to the Leibnitz rule and to
act on the generators via the usual formula
\be
\partial_r(b)={b\over r},\,\partial_r(h)={h\over r},\,\partial_r(c)={c\over r}.\label{deri} \ee
Recall that $r$ appears in the definition of a quantum
hyperboloid, and observe that this way of introducing the
derivative $\partial_r$ is compatible with the defining relations
of the algebras $\K_q[\R^3]$ and $\K_q[\R^4]$. Also, note that
formulae (\ref{deri}) are invariant w.r.t. a renormalization $r\to
a\,r,\,\,a\not=0$. So, a way of normalizing the Casimir  in the
definition of the q-hyperboloid does not matter.

Observe, the partial derivatives
$\partial_b,\,\partial_h,\,\partial_c$ that on the space $\SL$
can be identified with "bra" operators via the pairing
(\ref{pair}). Namely, we have
$$\partial_b=q\langle c,\,\, \partial_h={q^2\over 2_q} \langle h,\,\, \partial_c=q^3 \langle b.$$
Now, we look for coefficients $\mu$ and $\nu$ such that
 the following holds on $\SL$:
\be
\langle b=\mu {\B_q\over r^2}+\nu {b\over r}\partial_r,\quad \langle h=\mu {\H_q\over r^2}+\nu{h\over r}\partial_r,\quad
\langle c=\mu {\C_q \over r^2}+\nu {c\over r}\partial_r. \label{qraven}
\ee
\begin{proposition} These relations are valid with $\mu=q^{-4}$ and $\nu=q^{-2}$.
\end{proposition}
{\bf Proof} We have to check that by applying the operator equalities (\ref{qraven}) to any element
from $\SL$ we get a correct result.
For example, by applying the first equality to the element $b$ we get a correct relation
$$0=\langle b,b\rangle={q^2hB_q(b)-bH_q(b)\over q^4r^2}+{b\over rq^2} \partial_r(b)=0.$$
It is more difficult to check  three relations containing the
pairings $\langle b,c\rangle,\, \langle c,b\rangle$ and $\langle
h,h\rangle$. Let us check the first of them by leaving the others
to the reader:
$$q^{-3}=\langle b,c\rangle={q^2hB_q(c)-bH_q(c)\over q^4r^2}+
{b\over rq^2} \partial_r(c)={1\over q^4r^2}(q^2h({qh\over 2_q})-b(-c))+{1\over q^2r^2}bc=$$
$${h^2 \over q 2_qr^2}+
{bc \over q^4r^2}+{1\over q^2r^2}(cb+{1-q^2\over 1+q^2}h^2)={1\over q^3 r^2}(\qq bc+{h^2\over 2_q}+qcb)=q^{-3}.$$

Note, that the r.h.s. of the relations (\ref{qraven}) are well defined on the whole algebra $\K_q[\R^3]$ and
 extend the operators $\langle b,\, \langle h,\, \langle c$ to the algebra $\K_q[\R^3]$ via these relations.
Finally, we arrive to the following definition of the Laplace
operators in the algebras $\K_q[\R^3]$ and $\K_q[\Hr],$
respectively:
$$\De_{\K_q[\R^3]}=\qq \langle b \langle c+{1\over 2_q}\langle h \langle h +
q \langle c \langle b=q^{-5}\partial_c \partial_b+2_q q^{-4}\partial_h
\partial_h+q^{-3}\partial_b \partial_c, $$
$$\De_{\K_q[\Hr]}={1\over q^8 r^4}(\qq \B_q\C_q +{\H_q^2\over 2_q}+q\C_q\B_q).$$

Now, we define the braided Maxwell operators on these algebras as follows:
$$ \Mw_{\K_q[\R^3]}\left(\begin{array}{c}
\al\\
\beta\\
\gamma
\end{array}\right)=
\left(\begin{array}{c}
\De_{\K_q[\R^3]}(\al)\\
\De_{\K_q[\R^3]}(\beta)\\
\De_{\K_q[\R^3]}(\gamma)
\end{array}\right)-\left(\begin{array}{c}
\partial_b\\
\partial_h\\
\partial_c
\end{array}\right)(q^{-5}\partial_c,\,q^{-4}{2_q}\partial_h,\,q^{-3}\partial_b)\left(\begin{array}{c}
\al\\
\beta\\
\gamma
\end{array}\right).$$

$$\Mw_{\K_q[\Hr]}\left(\begin{array}{c}
\al\\
\beta\\
\gamma
\end{array}\right)=e\left(
\left(\begin{array}{c}
\De_{\K_q[\Hr]}(\al)\\
\De_{\K_q[\Hr]}(\beta)\\
\De_{\K_q[\Hr]}(\gamma)
\end{array}\right)-{1\over q^8 r^4}\left(\begin{array}{c}
\qq\C_q\\
{\H_q\over 2_q}\\
q\B_q
\end{array}\right)
(\B_q,\,\H_q,\,\C_q)\left(\begin{array}{c}
\al\\
\beta\\
\gamma
\end{array}\right)\right)$$
where $\al, \beta, \gamma$ in these formulas belong to the
algebras $\K_q[\R^3]$ and $\K_q[\Hr],$ respectively,
$e=1-{\overline e}$, and
$${\overline e}={1\over r^2}\left(\begin{array}{c}
\qq c\\
{h\over 2_q}\\
q b
\end{array}\right)
\left(\begin{array}{ccc}
b&h&c
\end{array}\right).
$$
Moreover, in the second formula we assume that the triples $(\al,
\beta,\gamma)^t$  belong to the projective module
$e\,\K_q[\Hr]^{\oplus 3}$ whereas in the first formula the columns
$(\al, \beta,\gamma)^t$ form the free module $\K_q[\R^3]^{\oplus
3}$.

Let us justify this definition.
The column $(\partial_b,\,\partial_h,\,\partial_c)^t$ in the first formula is universal (in the classical case
it corresponds to the de Rham operator).
The line  $(q^{-5}\partial_c,\,q^{-4}{2_q}\partial_h,\,q^{-3}\partial_b)$ is chosen  so that
$$\De_{\K_q[\R^3]}=(q^{-5}\partial_c,\,q^{-4}{2_q}\partial_h,\,q^{-3}\partial_b)\left(\begin{array}{c}
\partial_b\\
\partial_h\\
\partial_c
\end{array}\right)
$$

The second formula can be obtained from the first one by
disregarding the second summands in formulae (\ref{qraven}). Now,
we treat the module $e\,\K_q[\Hr]^{\oplus 3}$ as the space of
braided differential forms $\Om^1_q(\Hr)$
 (for other modules in question it can be
done in a similar way). The space $\Om^1_q(\Hr)$  treated as a
right $\K_q[\Hr]$ module is the quotient of the space of all
braided differential forms $(db)\, \al+(dh)\, \beta+(dc)\,\gamma$
over the forms $(\qq(db)\,c+{{(dh)\, h}\over 2_q}+q (dc)\,
b)\rho,\,\,\al, \beta, \, \gamma,\, \rho\in \K_q[\Hr]$. Note that
we have defined the element $\qq(db)\,c+{{(dh)\, h}\over 2_q}+q (dc)\, b$
by replacing the left factors in the Casimir ${\Cas}_{sl}$ by their "differentials"
without using either Leibnitz rule or any transposing "functions" and "differentials".

Similarly to the classical pattern, we have the following.

\begin{proposition} 1. The triples $(\partial_b\rho,\,\partial_h\rho,\,\partial_c\rho)^t$
belong to $\Ker(\Mw_{\K_q[\R^3]})$ provided the operator
$\De_{\K_q[\R^3]}$ commutes with
$\partial_b,\,\partial_h,\,\partial_c$.

2. The triples $(\qq\C_q\rho,\,{\H_q\over 2_q}\rho,\,q\B_q\rho)^t$ belong to $\Ker(\Mw_{\K_q[\Hr]})$ provided
 $$e\left(\left(\begin{array}{c}
\qq\C_q\\
{\H_q\over 2_q}\\
q\B_q
\end{array}\right)\De_{\K_q[\Hr]}-\De_{\K_q[\Hr]}   \left(\begin{array}{c}
\qq\C_q\\
{\H_q\over 2_q}\\
q\B_q
\end{array}\right) \right)=0$$
on the algebra $\K_q[\Hr]$.
\end{proposition}

Now, we pass to the q-Minkowski space algebra $\K_q[\R^4]$ and
 discuss a way to convert the Casimir element (\ref{casgl}) into an operator.
 Since ${\Cas}_{gl}={\Cas}_{sl}+{l^2\over 2_q}$ and a method
of assigning of an operator to the element ${\Cas}_{sl}$ is
presented above, we only have to define the operator $\partial_l$
on the algebra $\K_q[\R^4],$ keeping in mind, that
 we have $\partial_l={q^2\over 2_q}\langle l$ on the space
$\LL$. Since the element $l$ is central, it is possible to define
the extension of the derivative $\partial_l$ via the usual
Leibnitz rule. However, such a way is not compatible with the
first column of the system (\ref{sys2}). Nevertheless,
conjecturally  any element of the q-Minkowski space algebra
$\K_q[\R^4]$ can be written in a completely q-symmetrized form as
discussed in \cite{GPS2}, where this conjecture is proven for low
dimensional components. By assuming this conjecture to be true, we
define the derivative $\partial_l$ via the Leibnitz rule but only
on elements presented in such a "canonical" form.

\begin{remark} In general,  defining derivatives or vector fields  it is often convenient to do
this on  basis elements. This prevents us from checking that the
Leibnitz rule (if it is used in definition of  these operators) is
compatible with defining relations. This idea is close to that of
the paper \cite{G}, where  construction of the Koszul complexes
(similar to the de Rham complexes) uses  "R-symmetric" and
"R-skew-symmetric" algebras whose elements are realized in a
"(skew)symmetrized" form. It is also similar to the above
described construction  of the operators $B_q,\,H_q,\,C_q$.
\end{remark}

After having defined  this derivative we can present the Laplace
operator on the algebra in question in the following form
$$\De_{\K_q[\R^4]}=\De_{\K_q[\R^3]}+{1\over 2_q}\langle l \langle l=
\De_{K_q[\R^3]}+{2_q\over q^4}\partial_l \partial_l.$$
Also, define the Maxwell operator on this algebra
$$ \Mw_{\K_q[\R^4]}\left(\begin{array}{c}
\al\\
\beta\\
\gamma
\end{array}\right)=
\left(\begin{array}{c}
\De_{\K_q[\R^4]}(\al)\\
\De_{\K_q[\R^4]}(\beta)\\
\De_{\K_q[\R^4]}(\gamma)\\
\De_{\K_q[\R^4]}(\delta)
\end{array}\right)-\left(\begin{array}{c}
\partial_b\\
\partial_h\\
\partial_c\\
\partial_l
\end{array}\right)(q^{-5}\partial_c,\,q^{-4}{2_q}\partial_h,\,q^{-3}\partial_b,\,
{2_q\over q^4}\partial_l)\left(\begin{array}{c}
\al\\
\beta\\
\gamma\\
\delta
\end{array}\right).$$
Here $\al,\,\beta,\,\gamma,\delta\in \K_q[\R^4]$.

Now, we introduce an action of the QG $\uq$ on the ingredients of
the Maxwell operators in question so  that these operators become
$\uq$-invariant.

First, note that the QG acts on the operators $B_q,\,H_q,\,C_q$
and on $\langle b,\,\langle h,\,\langle c$ in the same way as on
the generators $b,\,h,\,c$. This means that the maps $b\to
B_q,...$ and $b\to \langle b,...$ are $\uq$  morphisms. So, we
restrict ourselves to the idempotent $\overline e$ and define the
above action such that this idempotent  becomes invariant. The
reader can easily extend our method to other ingredients of the
Maxwell operators.

Consider a representation of the algebra $U(sl(2)_q)$
$$\pi:b\to P^{-1}\, B_q\, P,\,\, h\to P^{-1}\, H_q\, P,\,\, c\to P^{-1}\, C_q\, P$$
where $P$ is an invertible numerical matrix: $P\in M_3(\K)$.
Define the action of the QG $\uq$ on the space $M_3(\K)$ according
to its action on the generators $b,\, h,\,c$:
$$K^\epsilon(P^{-1}\, B_q\, P)=q^{2\epsilon}P^{-1}\,
B_q\, P,...,Y(P^{-1}\, H_q\, P)=q^{\theta-1}2_q P^{-1}\, C_q\,
P,\,\,Y(P^{-1}\, C_q\, P)=0.$$ We extend this action to the whole
space $M_3(\K)$ by treating $M_3(\K)$ as the image  of the algebra
$U(sl(2)_q)$. By construction, the representation $\pi:U(sl(2)_q)
\to M_3(\K)$ is a $\uq$ morphism.

Consider the element
$$w^{-1}\pi_2({\Cas}_{sl})=
\qq bP^{-1}\left(\begin{array}{ccc}
0&0&0\\
-{ q\over 2_q}&0&0\\
0&q^2&0
\end{array}\right)P+{h\over 2_q}P^{-1}\left(\begin{array}{ccc}
q^2&0&0\\
0&q^2-1&0\\
0&0&-1\end{array}\right)P$$ $$
+qcP^{-1}\left(\begin{array}{ccc}
0&-1&0\\
0&0&{q\over 2_q}\\
0&0&0
\end{array}\right)P={1\over 2_q}
P^{-1}\left(\begin{array}{ccc}
q^2 h&-2_q q c &0\\
-b&(q^2-1)h&q^2c\\
0&2_q q b&-h\end{array}\right)P.$$ Here $\pi_2$ means that we
apply the representation $\pi$ to the second factors of the {\em
split} Casimir element, i.e. the Casimir element regarded as an
element of $\SL\ot\SL$.

Now, introduce the matrix  $L=(w^{-1}\pi_2({\Cas}_{sl}))^t$.
Note the the powers of this matrix are transposed to the "powers" of the matrix
$w^{-1}\pi_2({\Cas}_{sl})$ where these "powers" are defined via putting one matrix
$w^{-1}\pi_2({\Cas}_{sl})$ "inside" of the other one. This procedure is described
 in detail in \cite{GLS}.
It is easy to see that these "powers" are $\uq$-invariant. Also,
note that the matrix $L$ obeys  a CH identity (cf. \cite{GLS,
GS2}). Thus,  introducing a new action of the  QG $\uq$ on the
elements $M\in M_3(\K)$ as follows
$$X\triangleright M=(X(M^t))^t,\qquad \forall X\in \uq$$
where $X(M)$ is the above defined  action, we see that the matrix
$L$ and all its powers are $\uq$-invariant.

Now, we claim that the idempotent $\overline e$ can be presented as  a degree 2 polynomial
of the matrix $L$ for a special
choice of the matrix $P$. Therefore, this idempotent is also  invariant.
We leave finding this polynomial and the corresponding
 matrix $P$ to the reader.

\end{document}